\documentclass[mathcal,12pt]{article}
\usepackage{amsmath}
\usepackage{amsfonts}
\usepackage{amssymb}
\usepackage{latexsym}
\usepackage{euscript}
\usepackage{a4wide}
\usepackage{diagrams}
\usepackage{mathrsfs}

\def\Fil{\mathrm{Fil}}
\def\ModS{(Mod /S)}
\def\ModFIS{(Mod \ FI/S)}
\def\ModSp{(Mod /S_1)}
\def\ove{\overline{\ve}}
\def\M{\mathcal M}
\def\Lt{\mathcal L}
\def\B{\mathcal B}
\def\A{\mathcal A}
\def\G{\mathscr G}
\def\Ht{\mathscr H}
\def\D{\mathcal D}

\def\ve{\mathbf{e}}

\def\Kp{K^{+}}

\def\Eis{\mathfrak{I}}
\def\Eisv{\mathfrak{I}}

\def\Vbar{\overline{V}}
\def\mI{\mathfrak{m}}
\def\etale{\'{e}tale \ }
\def\etalestop{\'{e}tale}
\def\Bchi{B_{2,\omega^{k-2}}}
\def\Bka{B_{2k}}
\def\Bkb{B_{p+1-2k}}
\def\Bk{B_k}

\def\ra{\rightarrow}
\def\frak{\mathfrak}

\def\text{\textrm}

\def\GL{\mathrm{GL}}

\def\Qbar{\overline{\Q}}
\def\Gal{\mathrm{Gal}}
\def\F{\mathbf F}
\def\Fbar{\overline{\F}}

\def\rhobar{\overline{\rho}}
\def\varrhobar{\overline{\varrho}}
\def\rhouniv{\rho^{\mathrm{univ}}}

\def\qed{\hfill \square \ }

\def\Ot{\mathcal O}

\def\Q{\mathbf Q}
\def\Z{\mathbf Z}

\def\T{\mathbf T}
\def\twT{\widetilde{\T}}

\def\GL{\mathrm{GL}}

\def\Ext{\mathrm{Ext}}
\def\End{\mathrm{End}}

\def\Hom{\mathop{\mathrm{Hom}}\nolimits}

\def\Frob{\mathop{\mathrm{Frob}}\nolimits}
\def\Trace{\mathop{\mathrm{Trace}}\nolimits}
\def\Det{\mathop{\mathrm{Det}}\nolimits}

\newarrow{equals}{=}{=}{=}{=}{=}

\newtheorem{theorem}{Theorem}[section]
\newtheorem{lemma}[theorem]{Lemma}
\newtheorem{df}[theorem]{Definition}

\newtheorem{cor}[theorem]{Corollary}

\def\st{\star}

\def\XYZ{XY\kern-.17em{Z}}
\def\fppf{f\kern-.02em{p}p\kern-.07em{f}}
\def\Kf{K_{\kern-.1em{f}}}

\renewcommand\pmod{\mod}
\def\f1{f'}

\title{Eisenstein Deformation Rings}
\author{Frank Calegari\footnote{Supported in part by the American
Institute of Mathematics}}
\begin{document}
\maketitle
\abstract{We prove $R = \T$ theorems for certain reducible
residual Galois representations. We answer in the positive
a question of Gross and Lubin on whether certain Hecke
algebras $\T$ are discrete valuation rings. In order to
prove these results we determine (using the theory
of Breuil modules) when two finite flat group schemes
$\G$ and $\Ht$ of order $p$ over an arbitrarily
tamely ramified discrete valuation ring admit an extension not
killed by $p$.}

\section{Introduction}

In a previous paper \cite{CE}, M.~Emerton and the author studied 
modular deformation problems associated to
certain
\emph{reducible} representations. In
particular, for odd primes $p$
we considered   the totally split 
representation $\rhobar$ given by
$$\left(\begin{matrix} 1 & 0 \\ 0 & \chi \end{matrix}\right) \pmod p$$
where $\chi$ was the $p$-adic cyclotomic character.
It was proved in~$\cite{CE}$ that deformations
of $\rhobar$ finite flat at $p$ and satisfying a 
certain ``semistability'' condition at an auxiliary prime $N$
were modular of level $\Gamma_0(N)$, and the associated
universal deformation ring $R$ was isomorphic to
$\T_{\Eis}$, where $\T$ was the full Hecke algebra of level $\Gamma_0(N)$
and $\Eis$ was the $p$-Eisenstein ideal.
This enabled us to study the Eisenstein ideal by directly
studying deformations of $\rhobar$. 
In this paper, which can be seen as a sequel to
\cite{CE}, we study 
\emph{non-split} reducible representations
$\rhobar$ which are ramified only at $p$.
Under certain natural hypotheses, these representations are modular,
and
arise from cuspidal modular forms of
weight two and level $\Gamma_1(p)$ or $\Gamma_0(p^2)$.
We define certain deformation problems for $\rhobar$
such that the associated universal deformation ring
can be identified with the appropriate Hecke algebra
localized at the Eisenstein ideal, and 
use this to deduce properties of the Hecke algebras
in these cases.

\medskip

There are several important differences in the techniques of
this paper to those of~\cite{CE}.
In~\cite{CE}, the residual Eisenstein representation of level $\Gamma_0(N)$
is not minimal --- it has Serre conductor $1$. Thus one could play off
the minimal and non-minimal deformation problems using techniques of Wiles.
In this paper, $\rhobar$ is minimal of the appropriate level,
and the minimal deformation problem is not trivially $\Z_p$, as it was
in~\cite{CE}. The techniques used to prove modularity 
in this paper are quite different.
One ingredient is the following trivial observation. Suppose the following
are true:
\begin{itemize}
\item[(i)] $R \rightarrow \T$ is surjective,
\item[(ii)]
$\T \ne \T/p^n$ for any $n$,
\item[(iii)]  $R$ is a discrete valuation ring.
\end{itemize}
Then
$R = \T$.
By considering some very general properties of
residual representations we establish a criterion
that allows us to establish in many cases
that $R$ \emph{is}
a discrete valuation ring. For universal
deformation rings $R$ for which this criterion
does not apply, we construct another universal
deformation ring
$R'$ (corresponding to certain modular forms of
level $1$) such that we may apply our
criterion to  deduce that
$R'$ is a discrete valuation ring.
We then prove that
 $R/p = R'/p = \T/p$, and deduce from this
that $R = \T$.

\medskip

In our deformation problems we  consider finite flat group
schemes over bases $\Ot_K$ such that $e(K) > p - 1$, and thus we are forced
to utilize the theory of Breuil modules~\cite{Breuil}. In particular, we need
to consider 
finite flat group schemes $\G/\Ot_K$ that are not killed by $p$,
which leads to certain delicate computations with modules
over divided power rings.
As a consequence,
however, 
we prove an independently interesting result about certain group schemes
of order $p^2$ (see Theorem~\ref{theorem:oort1}).

\medskip

As in~\cite{CE}, the results of Skinner and Wiles~\cite{SW}
proving $R = \T$ theorems for reducible representations
do not apply, since our representations $\rhobar$ are
either locally split or are associated to non-ordinary
deformation problems.

\medskip

For the case of representations $\rhobar$ of level $\Gamma_1(p)$,
we are forced to make certain divisibility assumptions on
Bernoulli numbers. Probabilistically, these assumptions
should fail at most finitely often, but we have no proof of
this fact. It is clear that certain assumptions are required,
however. For example, without the assumption that the
the $\chi^{1-k}$-eigenspace of the class group of $\Q(\zeta_{p})$
is cyclic
it is not even clear
that the residual representations $\rhobar$ we consider are 
modular. These assumptions plays no role in the
$\Gamma_0(p^2)$ case, however.

\section{Results}

Let $p$ be an odd prime, let $K = \Q(\zeta_p)$,
and let $\Kp = \Q(\zeta_p)^{+}$ be the totally
real subfield of $K$.
Fix once and for all an embedding $\Qbar \rightarrow \Qbar_p$,
and let $K_p$ and $\Kp_p$ denote the respective
localizations
of $K$ and $\Kp$ inside $\Qbar_p$.
Let $2 < k < p -1$ be an even integer,  and
suppose that the $\chi^{1-k}$-eigenspace
of the class group of $K$ is non-trivial.
This is a well known consequence of Vandiver's
conjecture (see for example~\cite{Greenberg}).
Let $\chi$ be the cyclotomic character, and let
$\omega \equiv \chi \mod p$ be the Teichm\"{u}ller character.
Let $B_n$  denote the classical Bernoulli numbers, defined as follows:
$$ \frac{t e^t}{e^t - 1} = \sum_{n=0}^{\infty} \frac{ B_n t^n}{n!}.$$
If $\varphi$ is a character of $\Gal(\Qbar/\Q)$ of conductor $p$
then one may also define Leopoldt's generalized Bernoulli numbers $B_{n,\varphi}$
by the following generating series:
$$\sum_{k=0}^{p-1} \frac{\varphi(k) t e^{kt}}{e^{pt} - 1} = \frac{ B_{n,\varphi} t^n}{n!}.$$
For integral $n \ge 1$ one has the following congruence (Thm 2.3 of~\cite{Lang})
$$\frac{1}{n} B_{n,\omega^{k-n}} \equiv \frac{1}{k} B_k \mod p.$$

\medskip

Let $\rhobar$ be the unique non-split representation
$$\left(\begin{matrix} 1 & \st 
 \\ 0 & \chi^{k-1} \end{matrix}\right) \mod p,$$
that is unramified over $K$. It follows 
from~\cite{Ribet} that $\rhobar$ is modular
of weight $2$ and level $\Gamma_1(p)$.
Let $\Vbar$ be the two dimensional 
vector space on which $\rhobar$ acts.
We consider the following deformation problem
for $\rhobar$:
For a local artinian ring $A$ with residue field $\F_p$, let $\D(A)$ denote
the set of deformations $(\rho,V)$ of $(\rhobar,\Vbar)$
satisfying the following
properties:

\begin{enumerate}
\item The determinant of $\rho$ is $\chi \cdot \omega^{k-2}$.
\item $\rho$ is unramified outside $p$.
\item The representation $\rho_{\Kp_p}$ on $V$ is the generic fibre
of a finite flat group scheme $\G/\Ot_{\Kp_p}$.
\end{enumerate}
Since $e(\Kp_p) = \frac{1}{2} (p-1) < p-1$, the finite flat group scheme
$\G$ is determined up to isomorphism by $\rho$.
The following result is standard.
\begin{theorem} The functor $\D$ is $($pro$)$-representable
by a universal deformation ring $R$.
\end{theorem}
\begin{df} Fix $k$.
Let $\T$ be the cuspidal Hecke algebra of weight
$2$ and level $\Gamma_1(p)$
generated by $T_{\ell}$ for $\ell \ne p$. Then the $(p,k)$-Eisenstein
Ideal $\Eis$ is the maximal ideal of $\T$ containing
$T_{\ell} - 1 - \ell^{k-1}$ for all $\ell \ne p$. 
\end{df}
Since it will be clear from the context, we usually
refer to $\Eis$ as the Eisenstein ideal. Note that
we define $\Eis$ to be maximal, contrary to the usual
definition. Since we are only interested in $\T_{\Eis}$,
however, there should hopefully be no confusion.
We prove the following theorem.
\begin{theorem} \label{theorem:mainfirst}
Suppose that either $p \| \Bchi$ or $p \| \Bk$. Then
there is a natural isomorphism $R \simeq \T_{\Eis}$.
The ring $R$ is monogenic over $\Z_p$, and if
 $p \| \Bchi$, then $R$ is a discrete valuation ring.
\end{theorem}

\medskip

The other residual representations we consider in
this paper are wildly ramified, and arise from 
level $\Gamma_0(p^2)$.
Let $p$ be prime, and let $k < p-1$ be a positive integer
such that
$$k \ne 0,1,\frac{p-1}{2},\frac{p+1}{2}.$$
Let $k' < p-1$ be the unique positive integer
such that $k' + k \equiv \frac{1}{2}(p+1) \mod p-1$, and assume that
$p \nmid \Bka$ and $p \nmid \Bkb$.
Let $\rhobar$ be the unique non-split representation of the form
$$\left(\begin{matrix} \chi^k & \st 
 \\ 0 & \chi^{1-k} \end{matrix}\right) \mod p,$$
unramified away from $p$.  The Bernoulli
condition ensures that any non-split representation
is wildly ramified at $p$. The existence and uniqueness of $\rhobar$
is a simple exercise in class field theory. The representation
$\rhobar$ is modular of weight $2$ and level $\Gamma_0(p^2)$,
and if $2k > p+1$ actually occurs as a sub-representation of the
Jacobian $J_0(p^2)$ (see~\cite{Gross}). 
The Bernoulli number condition also ensures
that $\rhobar$ is
not a twist of a representation coming from
$\Gamma_1(p)$ (necessarily of the other residual representations we
are considering), and thus $\rhobar$ is ``genuinely'' of
level $p^2$.
Let $K/\Q_p$ denote a tamely ramified extension
of degree $p+1$.
Let $\Vbar$ be the two dimensional vector space on which $\rhobar$
acts.
We consider the following deformation problem
for $\rhobar$:
For a local artinian ring $A$ with residue field $\F_p$, let $\D(A)$ denote
the set of deformations $(V,\rho)$ of $(\Vbar,\rhobar)$ to $A$ satisfying the following
properties:
\begin{enumerate}
\item The determinant of $\rho$ is $\chi$.
\item $\rho$ is unramified outside $p$.
\item The representation $\rho|_{K}$ on $V$ is the generic fibre of a finite flat
group scheme $\G/\Ot_K$.
\end{enumerate}
Since $e(K) = p+1 \ge p-1$, finite flat group 
schemes are typically not determined by their generic fibre. It transpires,
however, that for the choice of $k$ above (in particular $2k
\not\equiv 0,2 \mod (p-1)$) that $\rhobar$ \emph{does} uniquely
determine a finite flat group scheme $\G/\Ot_K$. We
have the following:
\begin{theorem} The functor $\D$ is $($pro$)$-representable
by a universal deformation ring $R$.
\end{theorem}
\begin{df}
Let $\T$ be the cuspidal Hecke algebra of weight
$2$ and level $\Gamma_0(p^2)$. Then the $(p,k)$-Eisenstein
Ideal is the maximal ideal $\Eis$ of $\T$ containing
$T_{\ell} - \ell^k - {\ell}^{p-k}$ for all $\ell \ne p$.
\end{df}
We prove the following theorem.
\begin{theorem} \label{theorem:maingross}
There is a natural isomorphism $R \simeq \T_{\Eis}$.
The ring $R$ is a discrete valuation ring.
If $p \equiv 3 \mod 4$ and  $k = k' = (3p-1)/4$,
then $R \simeq \Z_p$.
\end{theorem}
This theorem was the main motivation for this paper.
It answers a question of Gross and Lubin, who asked
(\cite{Gross} p.310) whether $\T_{\Eis}$ was always
a discrete valuation ring. 

\medskip

Our last result is a consequence of finite flat
group scheme calculations required to prove Theorem~\ref{theorem:maingross},
although it is interesting in its own right.
First recall that (after choosing a uniformizer $\pi \in \Ot_K$) a
finite flat group scheme $\G/\Ot_K$ of order $p$ is specified
by its Oort--Tate parameters, a pair 
 $(r,a)$ with $r \in \Z$ satisfying
$0 \le r \le e = e(K/\Q)$, and $a \in \Ot_K/\frak{m}$.

\begin{theorem} \label{theorem:oort1}
 Let $K/\Q_p$ be a finite extension of ramification
degree $e$, where $(e,p) = 1$. Consider an exact sequence of finite
flat group schemes:
$$0 \rightarrow \Ht' \rightarrow \Ht \rightarrow \Ht'' \rightarrow 0$$
and suppose that $\Ht'$ and $\Ht''$ are finite flat group schemes of
order $p$. Then either
\begin{enumerate}
\item $\Ht$ is killed by $p$.
\item $\Ht$ is \'{e}tale or multiplicative, and $\Ht[p]$ is finite flat.
\item There exists a non-trivial morphism $\Ht'' \rightarrow \Ht'$
that is \emph{not} an isomorphism, and $\Ht[p]$ is not a finite flat
group scheme.
\end{enumerate}
Moreover, given a non-trivial morphism $\Ht'' \rightarrow \Ht'$
 there exists an extension
$\Ht \in \Ext(\Ht'',\Ht')$ such that $\Ht \ne \Ht[p]$ if and only
if the Oort--Tate parameters $(s,b)$ of $\Ht'$ and
$(r,a)$ of $\Ht''$
 satisfy either of the following inequalities:
$r \ge ps$ or $(e-s) \ge p(e-r)$.
\end{theorem}

\medskip

I would like to thank Brian Conrad and Christophe Breuil
for answering some technical
questions about finite flat group schemes
and Breuil modules respectively.
I would also like to thank Matthew Emerton for
our frequent conversations, during which several
of the ideas of this paper had their genesis.

\section{Generalities on Eisenstein representations}
\label{section:chen}

The main reference for this section is the
paper of Bella\"{\i}che and Chenevier~\cite{Chenevier}.
In this section we record some general 
remarks about residually reducible representations. 
Let $(A,\mI,k)$ be a local $p$-adically complete
ring.
Given a  representation $\rho: \Gal(\Qbar/\Q)
\rightarrow \GL_2(A)$ such that 
$$\rhobar: \Gal(\Qbar/\Q) \rightarrow \GL_2(A/\mI) = \GL_2(k)$$
is reducible and unramified outside $p$ and $\infty$
 we shall derive a sufficient criterion for  
$A$ to be a discrete valuation ring.
The results of this section should not be considered original,
and follow almost directly from the arguments of~\cite{Chenevier}.
The spirit of these arguments is also very similar to the
work of
Papier and Ribet~\cite{Papier}.
We shall use the notation of~\cite{Chenevier}, however.

\medskip

Let $G$ be a quotient of $\Gal(\Qbar/\Q)$ such that
$\rho$ factors through $G$. We may  consider
$\rho$ as a representation of $G$ into $\GL_2(A)$.
Let $T: G \rightarrow A$ denote the composite of $\rho$
with the trace map. Suppose that the semi-simplification
$(\rhobar)^{ss}$ is given by $\chi_1 \oplus \chi_2$.
We shall assume that $\chi_1 \ne \chi_2$.
Fix $s \in G$ such that $\chi_1(s) \ne \chi_2(s)$. 
The characteristic
polynomial of $\rho(s)$ has two distinct roots modulo $\mI$, and
thus by Hensel's Lemma has roots $\lambda_1$ and $\lambda_2$ with
$\lambda_i \equiv \chi_i(s) \mod \mI$. Choose a basis
of the representation $\rho$ such that
$\rho(s) e_i = \lambda_i e_i$. Let $a$, $b$, $c$, $d$ be
the matrix entries of $\rho$ with respect to this basis,
and let $B$ and $C$ be the $A$-ideals generated by
$b(g)$ and $c(g)$ respectively, for $g \in G$.
Let $I \subset A$ be a proper ideal such that
$T \mod I$ can be written as the sum of two characters
$\psi_1$, $\psi_2$ such that $\psi_i \mod \mI = \chi_i$.

\begin{lemma} For all $g$, $g' \in G$, 
$a(g) - \psi_1(g) \in I$, $b(g) - \psi_2(g) \in I$, and
$b(g)c(g') \in I$.
\end{lemma}

\begin{Proof} This is Lemme 1 of~\cite{Chenevier}.
\end{Proof}

\begin{lemma}  There is an injection of $A$-modules
$$\Hom_A(B/IB,A/I)  \rightarrow  \Ext^1_{(A/I)[G]}(\psi_2,\psi_1).$$
\end{lemma}

\begin{df} The ideal of reducibility of $A$ is the
largest ideal of $A$ such that $T \mod I$ is the
sum of two characters. There is an equality $I = BC$.
\end{df}

\begin{lemma} Suppose that $A$ is noetherian,
that  the ideal of
reducibility is maximal, and that
$$\dim_k \Ext^1_{k[G]}(\chi_2,\chi_1) = \dim_k \Ext^1_{k[G]}(\chi_2,\chi_1)
= 1.$$
 Then the maximal ideal $\mI$ of $A$ is
principal. If moreover $A$
admits a surjective map 
to a ring of characteristic zero,
then $A$ is a discrete valuation ring.
\label{lemma:3point4}
\end{lemma}

\begin{Proof} One has $\dim_k B \otimes_A k \le 1$.
Thus by Nakayama's Lemma, $B$ is a cyclic $A$ module,
and hence principal. A similar argument applies to $C$,
and thus $\mI = I = BC$ is principal.
Let $\mI = (\pi)$.
By Krull's Intersection Theorem each element of $A$ is of the
form $u \pi^k$ for some unit $u \in A$. Since $A$
admits a surjective  map to a ring of characteristic zero
$\pi$ is not nilpotent. Thus $A$ is a discrete valuation
ring.
\end{Proof}

\subsection{The General Strategy}

Let us now explain the general strategy of this paper.
In both cases we are considering a reducible representation
$\rhobar$, and a suitable universal deformation 
$$\rho^{univ}: \Gal(\Qbar/\Q) \rightarrow \GL_2(R)$$
unramified outside $p$ and $\infty$.
If $\Q_{\{p,\infty\}}$ denotes the maximal extension of $\Q$
unramified outside $p$ and $\infty$ then $\rhouniv$ 
\emph{a fortiori} factors through $\Gal(\Q_{\{p,\infty\}}/\Q)$.
Moreover, in either case
$\chi_1/\chi_2$ is some non-trivial power of the cyclotomic
character, and so  of the form $\chi^i$ for some $i \ne 0$.
  Our assumptions regarding Vandiver's
conjecture (in the case of $\Gamma_1(p)$, and automatically
in the case of $\Gamma_0(p^2)$)
imply that $\Ext^1_{\F_p[G]}(1,\chi^i)$ and $\Ext^1_{\F_p[G]}(1,\chi^{-i})$ are 
one dimensional, where
$G = \Gal(\Q_{\{p,\infty\}}/\Q)$. The universal deformation
rings are topologically finitely generated over $\Z_p$ and thus
noetherian.  Thus if $I$ is the ideal of
reducibility of $R$, then $I$ is principal.
Providing $R$ admits a surjection
to $\T$, we infer that $R$ is a discrete valuation
ring whenever the ideal of reducibility is maximal. 
The following lemma is trivial:

\begin{lemma} The ideal $I$ of $R$ is maximal 
if and only if there does not exist a surjection
$R/I \rightarrow \F_p[x]/x^2$ or $R/I \rightarrow \Z/p^2 \Z$.
\end{lemma}

In view of the description of $R$ as a universal deformation
ring, it therefore suffices to show that $\rhobar$ does not admit
any non-trivial deformations to $\GL_2(\F_p[x]/x^2)$ or $\GL_2(\Z/p^2 \Z)$ 
that are upper triangular. For the representations we consider 
of level $\Gamma_0(p^2)$,
it turns out that there are never any such deformations.
For level $\Gamma_1(p)$, however,
their may exist upper triangular
deformations to $\Z/p^2 \Z$.
This happens whenever
$p^2 | \Bchi$ (I do not know any example where
this happens, although it is conjectured to 
happen infinitely often).
To deal with this possibility, we switch to  another deformation
ring $R'$ corresponding to deformations of $\rhobar$ that 
arise
from modular forms of weight $k$ and level
$\Gamma_0(1)$. The ring $R'$ is a discrete
valuation ring whenever $p^2 \nmid \Bk$ by the same proof as
in Lemma~\ref{lemma:3point4}. By our assumptions
on Bernoulli number divisibility 
 this is always the case.
Thus the failure of
$R$ to be a discrete valuation ring
forces $R'$ to be a discrete valuation ring.
In this situation we find that
$R' \simeq \T'$, where
$\T'$ is the cuspidal Hecke algebra of weight $k$ and level $\Gamma_0(1)$
localized at the Eisenstein ideal. One knows, however, that
$\T'/p  \simeq \T/p$. Moreover, by purely local considerations it
follows that $R'/p \simeq R/p$. From these facts (along with
the observation that $\T$ is torsion free) we may conclude that
$R = \T$.

\medskip

One of the  main technical difficulties of the paper is 
determining the upper triangular deformations of $\rhobar$
to $\F_p[x]/x^2$ and $\Z/p^2 \Z$.
Note that it is \emph{not} always the case that the 
``Eisenstein Ideal'' as defined by Mazur (and others)
is the ideal of reducibility. Indeed, for $J_0(N)$ this
is never the case. For example, for $N=11$, the
Hecke algebra $\T \simeq \Z$ and the Eisenstein ideal
$\Eis \simeq (5)$. However, one easily finds that
the ideal of reducibility is $(25)$. This was noted
by Serre, and in the optic of the Eisenstein ideal was
pointed out by Mazur  (see for example the discussion
in~\cite{eisenstein}, Prop 18.9, pp.138--139). 

\medskip 

As a computational observation, it is typically the case that
$R = \T_{\Eis} = \Z_p$. This is not always true, however. 
For example,
when $p=547$ and $k=486$,
and $\rhobar$ is the residual representation
of level $\Gamma_1(547)$, then
using William Stein's Modular Forms Database one finds
that
$$\T_{\Eis}/p \simeq \F_p[x]/x^2.$$
Similarly, although I know of no  examples,
there is no reason why $p^2$ cannot divide
$\Bchi$. Note, however, that
if both conditions occur simultaneously, then $R$
\emph{cannot} be a discrete valuation ring. This follows
from the fact that the ideal of reducibility $I = (a)$ is
principal, and
the only discrete valuation ring $R$
that admits a surjection $R/a \rightarrow \Z/p^n$ for
some $n \ge 2$ is $\Z_p$.

\section{Deformations of level $\Gamma_1(p)$ and $\Gamma_0(p^2)$.}

\subsection{Eisenstein Deformations at level $\Gamma_1(p)$.}
Let $2 \le k < p-1$.
The residual representations considered in this section
were studied by Ribet~\cite{Ribet}, who proved that
whenever $p | \Bk$ for even $k$, the $\chi^{1-k}$-eigenspace of
the class group of $K = \Q(\zeta_p)$ is non-trivial
(The converse of this statement is a more classical theorem of
Herbrand). Moreover, given such a $p$ and $k$,
there exists a non-split representation
$$\left(\begin{matrix} 1 & \st 
 \\ 0 & \chi^{k-1} \end{matrix}\right) \mod p,$$
that is unramified over $K$. The arguments of
Ribet may be summarized as follows.
Since $B_2 = 1/6$, $k \ne 2$ and thus
$\chi^{k-1} \ne \chi$. Inside some variety
$J$ isogenous to $J_1(p)/J_0(p)$ one may
find a non-split representation $\rhobar$ of the shape above.
Now $J$ acquires
 everywhere good reduction over the
totally real subfield $K^{+}$ of $K$. The
results of Raynaud~\cite{Raynaud} imply that group schemes over
a base of ramification $e < p-1$ are determined
by their generic fibre. Thus the representation
$\rhobar$ over $\Kp_p$ is seen to arise
from a finite flat group scheme over $\Ot_{\Kp_p}$
that is an extension of a local
group scheme by an \'{e}tale group scheme. The connected--\'{e}tale
sequence therefore splits this extension of group schemes,
and thus $\rhobar$ is locally split at $p$ over $K$.
This implies that $\rhobar$ is the representation
considered above.
Let $\Vbar$ be the two dimensional representation
corresponding to $\rhobar$.
For a local artinian ring $A$ with residue field $\F_p$, let $\D(A)$ denote
the set of deformations $(V,\rho)$ of $(\Vbar,\rhobar)$ 
to $A$ satisfying the following
properties:
\begin{enumerate}
\item The determinant of $\rho$ is $\chi \cdot \omega^{k-2}$.
\item $\rho$ is unramified outside $p$.
\item The representation $\rho|_{\Kp_p}$ on $V$  arises from a finite flat group scheme over
$\Ot_{\Kp_p}$.
\end{enumerate}
\begin{theorem} The functor $\D$ is $($pro$)$-representable
by a universal deformation ring $R$.
\end{theorem}

\begin{Proof} 
This result can now be
considered relatively standard (see for 
example~\cite{CDT} and~\cite{BCDT}).
Let us make a few remarks, however.
First note
that $\End_{\Gal(\Qbar/\Q)}(\Vbar) = \F_p$, and thus
the ``unadorned'' universal deformation ring exists.
Let us make precise what is meant  by saying that $\rho_{\Kp_p}$
arises from a finite flat group scheme over $\Ot_{\Kp_p}$.
Essentially it stipulates the existence of a group scheme
$\G/\Ot_{\Kp_p}$ such that the induced representation
of $\Gal(\Qbar_p/\Kp_p)$ on the generic fibre gives rise to
$\rho|_{\Kp_p}$. Since $\rho$ is defined over $\Q_p$,
this implies that the generic fibre of $\G$ also descends
to $\Q$. Thus we automatically obtain a pair 
$(\G,\phi)$, where $\G$ is a finite flat group scheme over
$\Kp$, and $\phi$ is an action of $\Gal(K^{+}/\Q)$ on
the generic fibre of $\G$ (which extends to an action
on $\G$). One calls such a pair $(\G,\phi)$ a group
scheme with \emph{generic fibre descent data}. 
Note that since $e < p-1$, the group schemes $\G$
are uniquely determined by the deformations $\rho$.
$\qed$
\end{Proof}

\begin{lemma} There is no non-trivial element of $\D(F_p[x]/x^2)$ such
that $\rho$ is upper triangular. The ring $R$ is generated
as a $\Z_p$-algebra by traces.
\label{lemma:upper}
\end{lemma}

\begin{Proof}
Suppose that $\rho: \Gal(\Qbar/\Q) \rightarrow \GL_2(\F_p[x]/x^2)$
is upper triangular, and let $\G/\Ot_{\Kp_p}$ be the associated
finite flat group scheme. Let $\psi$ be the character corresponding
to the upper left hand corner of the representation. Then the
Galois representation $\psi|_{\Kp_p}$ gives rise to a sub-representation
of the generic fibre of $\G$, and thus to a finite flat subgroup scheme
$\Ht$ of $\G$. The generic fibre of $\Ht$ has a filtration
by constant Galois modules.
Since $e = \frac{1}{2}(p-1) < p-1$, $\Ht$ is therefore
an extension of constant group schemes. 
Thus $\Ht$ is an extension of \etale group schemes, and 
therefore $\Ht$ is \etalestop.
Thus $\psi$ considered as a character of $\Gal(\Qbar/\Q)$ is unramified at $p$
and thus unramified everywhere. By simple class field theory
it follows that $\psi$ is trivial.
In particular, $\rho$ must have the shape:
$$\left(\begin{matrix} 1 & \st \\ 0 & \chi^{k-1} \end{matrix}\right) \in \GL_2(\F_p[x]/x^2).$$
As in Ribet~\cite{Ribet}, the connected--\etale sequence implies that $\G$ splits over
$\Ot_{\Kp_p}$, and thus that $\rho$ is unramified over $\Kp$.
If $\rho$ defines a non-trivial representation to $\GL_2(\F_p[x]/x^2)$, we see
that its kernel must cut out a $(\Z/p\Z)^2$ unramified extension of $\Q(\zeta_p)$.
Since this contradicts our assumptions on $\rhobar$, the result follows.
To show that $R$ is generated by traces, it suffices to show that
any non-trivial deformation of $\rhobar$ to $\GL_2(\F_p[x]/x^2)$ is
generated by traces. This follows in
a standard way from the fact that (by Nakayama's Lemma)
$R$ is generated as a $\Z_p$-algebra by the generators of
$\mI_R/(\mI^2_R,p)$. 
Let $\rho$ be a deformation of $\rhobar$  to $\GL_2(\F_p[x]/x^2)$
than cannot be written as an upper triangular representation.
Write
the matrix entries of $\rho$ as functions $a$, $b$, $xc$ and $d$ of
$\Gal(\Qbar/\Q)$.  Then if $\sigma \in \Gal(\Qbar/\Q(\zeta_p))$, then
$\Det(\rho(\sigma)) - \Trace(\rho(\sigma)) = x  b(\sigma) c(\sigma)$.
Since $c$ is non-trivial (by assumption), the Cebotarev
density theorem implies there exists a $\sigma$
such that
that $b(\sigma)c(\sigma) \ne 0$.
Since $\Det(\rho(\sigma))  = 1$, it
follows that the traces of $\rho$ generate $\F_p[x]/x^2$.
$\qed$
\end{Proof}

\begin{df}
Let $\T$ be the cuspidal Hecke algebra of weight
$2$ and level $\Gamma_1(p)$. Then the $(p,k)$-Eisenstein
Ideal $\Eis$ is the maximal ideal of $\T$ containing
$T_{\ell} - 1 - {\ell}^{k-1}$ for all $\ell \ne p$.
\end{df}

\begin{lemma} There exists a surjective map $R \rightarrow \T_{\Eis}$.
\label{lemma:surj}
\end{lemma}

\begin{Proof} If $\twT_{\Eis}$ denotes the
normalization of $\T_{\Eis}$, then we may write
$\twT_{\Eis} = \prod_{i=1}^d \Ot_i$, where each $\Ot_i$ is a discrete
valuation ring finite over $\Z_p$. The rings $\Ot_i$ are in bijection
with normalized newforms $f$ of level $\Gamma_1(p)$ such
that if $\rho_{f}$ is an integral $p$-adic Galois representation
associated to $f$, then $(\rhobar_f)^{ss}  = (\rhobar)^{ss}$.
Arguing as in Ribet~\cite{Ribet},
For each form $f$, there exists a lattice such that the
reduction $\rhobar_f$ is a non-split extension of
$\chi^{k-1}$ by $1$. Since $f$ has level $\Gamma_1(p)$
this implies that $\rhobar_f$  splits over $K_p$ and thus is
unramified after restriction to $K$. By assumption
(on the cyclicity of the $\chi^{1-k}$-eigenspace of the class group)
this uniquely determines $\rhobar_f$, and thus $\rhobar_f = \rhobar$.
It follows that $\rho_f$ is a deformation of $\rhobar$, 
and thus
there exists a map $R \rightarrow \Ot_i$. In particular, we obtain
a map $R \rightarrow \prod_{i=1}^d \Ot_i = \twT_{\Eis}$. Since
$R$ is generated by traces, the image $R$ is also generated
by traces. The image of $\Trace(\rho^{univ}(\Frob_{\ell}))$
 for $\ell \ne p$ is $T_{\ell} \in \T$.
Since Frobenius elements are dense in $\Gal(\Qbar/\Q)$,
the image of $R$ is exactly  $\T_{\Eis}$.
$\qed$
\end{Proof}

\begin{lemma} \label{lemma:lower} If $p^2 \nmid \Bchi$ then there are no
deformations of $\rhobar$  in $\D(\Z/p^2 \Z)$ that are upper triangular.
\end{lemma} 

\begin{Proof} As in Lemma~\ref{lemma:upper}, the character $\psi$ corresponding
to the upper left hand corner must be trivial,
and thus $\rho$ is of the form
$$\left(\begin{matrix} 1 & \st \\ 0 & \chi \omega^{k-2} \end{matrix}\right) \in \GL_2(\Z/p^2\Z).$$
where $\rho|_{\Q_p(\zeta_{p^2})}$ is totally split. Thus this defines a degree $p^2$ unramified
extension of $\Q(\zeta_{p^2})$, which implies the divisibility of
Bernoulli numbers. $\qed$
\end{Proof}

\begin{cor} Suppose that   $p^2 \nmid \Bchi$. Then
$R$ is a discrete valuation ring and $R = \T$.
\end{cor}
\begin{Proof}  By Lemmas~\ref{lemma:upper} and~\ref{lemma:lower}, we conclude
that the ideal of reducibility $I$ of $R$ is maximal. Since $R \rightarrow \T$,
$p$ is not nilpotent in $R$, and thus $R$ is a discrete valuation ring.
That $R \simeq \T$ is then obvious, since $\T$ is non-trivial and has
characteristic zero.
$\qed$
\end{Proof}

\medskip

Thus to complete the proof of Theorem~\ref{theorem:mainfirst} 
we are now left to consider the case that $p^2 | \Bchi$.
We may therefore assume that $p \| \Bk$.
Note that for any $p$, the na\"{\i}ve probability that
there exists a $2 \le  k < p-1$ such that $p^2 \| \Bchi$
is approximately $1/p$. The further condition that
$p^2 \| \Bk$ decreases this probability to $1/p^2$. Thus
one might suppose that the divisibilities $p^2 | \Bk$ and
$p^2 | \Bchi$ occur at most finitely often for all $p$.
I do not now of any examples in which either condition
is satisfied. 

\medskip

Let us now assume that $p \| \Bk$.
Instead of proving the modularity of $R$
directly, we shall switch to
another deformation problem. We note, firstly, that the
representation $\rhobar$ is modular of level $1$ and weight $k$.
Certainly there exists a non-split representation $\rhobar'$
of that level with $(\rhobar')^{ss} = (\rhobar)^{ss}$. On the other hand,
$\rhobar'$ is crystalline  (in the sense of Fontaine--Laffaille~\cite{FL})
) with Hodge--Tate weights $[0,k-1]$, which must necessarily be split
locally over $K_p$. By our assumptions such a representation is
unique, so it must  equal $\rhobar$.
 We define the following deformation problem $\D'$.

\begin{enumerate}
\item The determinant of $\rho$ is $\chi^{k-1}$.
\item $\rho$ is unramified outside $p$.
\item The representation $\rho$ is ordinary at $p$, i.e.
there exists an exact sequence:
$$0 \rightarrow V' \rightarrow V \rightarrow V'' \rightarrow 0.$$ 
where $V'$ and $V''$ are free $A$-modules of rank $1$,
and $\Gal(\Qbar_p/\Q_p)$
acts on $V''$ via an unramified character.
\end{enumerate}

\begin{theorem} The functor $\D'$ is representable by a universal
deformation ring $R'$. Moreover, if $A$ is an artinian ring killed
by $p$, then $\D(A) \subseteq \D'(A)$, and thus there is a surjection
$R'/p \rightarrow R/p$.
\label{theorem:4point7}
\end{theorem}

\begin{Proof}  The existence of $R'$ is standard.
One could also try and define $\D'$ to be deformations that
are crystalline with Hodge--Tate weights $[0,k-1]$ and
presumably this would be define an equivalent functor.
 Let $(\rho,V)$ be a deformation in $\D(A)$.
Consider the connected--\etale sequence attached
to the finite flat group scheme associated to $\rho$.
On generic fibres, it induces an exact sequence
$$0 \rightarrow V' \rightarrow V \rightarrow V'' \rightarrow 0.$$
Since $\Vbar''$ and $\Vbar'$ are one dimensional it follows from
Nakayama's Lemma that $V''$ and $V'$ are cyclic. By a counting
argument it follows that $V'$ and $V''$ are free. 
Since $e = \frac{1}{2}(p-1) < p$, this splitting
descends to $\Q_p$, and it follows that $V$ is ordinary.
The existence of a
surjection $R'/p \rightarrow R/p$ follows by Yoneda's Lemma.
$\qed$
\end{Proof}

Note that this argument also implies that
deformations   $\rho \in \D(A)$ are in general ordinary.
 If $\rho^{mod}$ is
the representation
$$\rho^{mod}: \Gal(\Qbar/\Q)\rightarrow
\GL_2(\T_{\Eis})$$
 it follows that
  $T_p$ is given by the action of Frobenius
on the unramified quotient. Thus $T_p \in \T_{\Eis}$.

\begin{lemma} if $p \|  \Bk$ then $R'$ is a discrete
valuation ring.
\end{lemma}

\begin{Proof} Let $I'$ be the ideal of reducibility of $R'$.
The case of upper triangular deformations to
$\GL_2(\F_p[x]/x^2)$ is essentially the same argument as
for $R$, except that now the splitting of Galois modules
over $K_p$ comes from the ordinary hypothesis rather
than the connected--\etale sequence.
Consider a reducible deformation $\rho \in \D'(\Z/p^2 \Z)$.
Let $\psi$ be the character corresponding to the upper left hand 
corner of the representation. Then $\psi|_{\Q_p}$ is an extension
of a trivial representation by a trivial representation. If $\psi$ is
ramified at $p$ then $\rho$ can certainly not be ordinary, so
$\psi$ is unramified at $p$ and thus trivial. 
Hence $\rho$ is of the form
$$\left(\begin{matrix} 1 & \st \\ 0 & \chi^{k-1} \end{matrix}\right) 
\in \GL_2(\Z/p^2 \Z).$$
Moreover the ordinary hypothesis implies that the representation
must split locally over $K_p$ 
The kernel of $\rho$ defines
a degree $p^2$ unramified $\chi^{1-k}$ extension of $\Q(\zeta_{p^2})$, and
in particular implies that $p^2 | \Bk$. This contradicts 
our assumption.
 Thus $I$ is maximal in $R'$, and thus
$R'$ is a discrete valuation ring. $\qed$
\end{Proof}

\medskip

Let  $\T'$ be the cuspidal
Hecke algebra of level one and weight $k$, and let
$\Eisv$ be the Eisenstein ideal. Since $k < p -1$, the cuspidal
Eisenstein deformations are ordinary and in the usual way we
obtain a surjection $R' \rightarrow \T'_{\Eisv}$, which must be
an isomorphism. Note that this is expected since ordinary
representations of weight $k > 2$ are automatically crystalline.
There is a standard identification $\T'_{\Eisv}/p \simeq
\T_{\Eis}/p$ which follows from the
identification $S_2(\Gamma_1(p),\omega^{k-2},\Z/p\Z)
= S_k(\Gamma_0(1),\Z/p\Z)$. Thus $R'/p = \T'_{\Eisv}/p = \T_{\Eis}/p
= \F_p[x]/x^e$ for some $e$ (recall that $R'$ is a discrete
valuation ring, so monogenic over $\Z_p$). Since
$R'/p$ surjects onto $R/p$ by Theorem~\ref{theorem:4point7} and $R/p$ surjects onto $\T_{\Eis}/p$,
it follows that $R/p = \F_p[x]/x^e$ also.
Thus there exists a diagram:
$$\begin{diagram} \Z_p[[x]] &  \rOnto & R & \rOnto & R/p \simeq \F_p[x]/x^e\\
   &  & \dOnto & & \dequals \\
  &  & \T_{\Eis} & \rOnto & \T_{\Eis}/p \simeq \F_p[x]/x^e. \end{diagram}$$
Since $\T_{\Eis}$ is a monogenic and torsion free, it must
be isomorphic to $\Z_p[[x]]/f$ for some polynomial $f(x) \equiv x^e \mod p$
of degree $\le e-1$. 
Thus $R \simeq \Z[[x]]/I$ with $I = Jf$.
 Since $R/p \simeq 
\F_p[x]/x^e$, it follows that the image of the ideal $J$
in $\F_p[x]/x^e$ contains $1$. Thus $J$ contains $1$, and $R \simeq 
\T_{\Eis}$. This completes the proof 
of Theorem~\ref{theorem:mainfirst}.

\medskip

It was observed by William Stein~\cite{calstein} that although
$p$ can divide the discriminant of the Hecke
algebra of weight $2$ and level $\Gamma_0(p)$ 
(for example when $p=389$), it never appears to divide
the index. Equivalently, if $\T_{\frak{m}}$ is a
localization of the Hecke algebra $\T$ of weight two
and level $\Gamma_0(p)$, then $\T_{\frak{m}}$ is a discrete
valuation ring.
Computations of Stein also suggest this conjecture
may be true at level $\Gamma_1(p)$. Although we do not
prove this conjecture in the Eisenstein case, there 
are some interesting connections that arise. 
If $E$ is an elliptic curve of conductor $p$
then the associated
Hecke algebra is an integral domain if and only if $p$
does not divide the modular degree. As Stein notes,
a result of Flach implies that
the modular degree annihilates a certain Selmer
group~\cite{Flach}, which can in turn  be considered a form
of generalized class group. Thus the conjecture that
$\T_{\frak{m}}$ is an integral domain translates
into the conjecture that $p$ does not divide the
order of a certain ``class group'', and so resembles
the statement of Vandiver's conjecture. 
For Eisenstein representations of level $\Gamma_1(p)$,
we see that the same question is intimately related
 to the actual
Vandiver's conjecture.

\subsection{Eisenstein Deformations at level $\Gamma_0(p^2)$.}
\label{section:gross}

Let $p$ be prime, and let $k < p-1$ be a positive integer
such that
$$k \ne 0,1, \frac{p-1}{2}, \frac{p+1}{2}.$$
Let $k' < p-1$ be the positive integer such
that $k + k' = \frac{1}{2}(p+1) \mod p-1$, and assume that
$p \nmid \Bka$ and $p \nmid \Bkb$.
Under these conditions, there is a unique
representation $\rhobar$  of the form
$$\left(\begin{matrix} \chi^k & \st 
 \\ 0 & \chi^{1-k} \end{matrix}\right) \mod p,$$
which is wildly ramified at $p$ and unramified outside $p$.
It follows from~\cite{Gross} that $\rhobar$
is modular of weight $2$ and level $\Gamma_0(p^2)$.
Let $K/\Q_p$ denote a tamely ramified extension
of degree $p+1$.
We shall study deformations of $\rhobar$ that
arise from finite flat group schemes over $K$.
Since $e(K) = p+1 > p-1$, however, a Galois
module that arises from a finite flat group scheme
does not necessarily determine a \emph{unique}
finite flat group scheme.
We prove, however, the following result (see
Section~\ref{section:breuil}, Lemma~\ref{lemma:theoremX}):
\begin{theorem} Let $k$ be as above,
and let  $V$ be the representation
corresponding to $\rhobar$.
Then there exists a unique finite
flat group scheme $\G/\Ot_K$ with generic fibre
descent data to $\Q_p$ isomorphic to $V$.
\end{theorem}

To prove this theorem, we need to study
the category of finite flat group schemes
over $K$. An explicit theory of
finite flat group schemes over discrete
valuation rings of arbitrary ramification
was constructed by Breuil~\cite{Breuil}.

\medskip 

For a local artinian ring $A$ with residue field $\F_p$, 
let $\D(A)$ denote
the set of deformations $(V,\rho)$ of $(\Vbar,\rhobar)$ 
to $A$ satisfying the following
properties:
\begin{enumerate}
\item The determinant of $\rho$ is $\chi$.
\item $\rho$ is unramified outside $p$.
\item The representation $\rho|_{K}$ on $V$ 
is the generic fibre of a finite flat
group scheme $\G/\Ot_K$.
\end{enumerate}

\medskip

Note that finite flat group schemes over
$\Ot_K$ do not form an abelian  category. However,
throughout this paper we implicitly use the following lemma:
\begin{lemma} Let
$$0 \rightarrow H' \rightarrow H \rightarrow H'' \rightarrow 0$$
be a short exact sequence of finite Galois modules such that
$H$ is the generic fibre of a finite flat group scheme
$\Ht/\Ot_K$. Then there exist (unique) finite flat group schemes
$\Ht'$, $\Ht''/\Ot_K$ which fit into a short exact sequence
$$0 \rightarrow \Ht' \rightarrow \Ht \rightarrow \Ht'' \rightarrow 0,$$
and such that taking generic fibres in this sequence yields the
exact sequence of Galois modules above.
\label{lemma:conrad}
\end{lemma} 
This lemma is proved in~\cite{Conrad} \S$1.1$. The group scheme
$\Ht'$ will be the scheme theoretic closure of $H'$ inside $\Ht$.
Note however that even if the map $H \rightarrow H''$
is ``multiplication by $p$'', this does not identify
$\Ht'$ with $\Ht[p]$. Indeed, it might be the case
that $\Ht[p]$ is not even a finite flat group scheme.
However, in view of Lemma~\ref{lemma:conrad} one certainly
has the following: 

\begin{theorem} The representation
$\D$ is representable by a universal deformation ring $R$.
\end{theorem}

\begin{Proof}
This follows in the standard way. $\qed$
\end{Proof}

\begin{lemma} There is no non-trivial element of $\D(\F_p[x]/x^2)$ such
that $\rho$ is upper triangular. The ring $R$ is generated
by traces, and moreover there exists a surjective map
$R \rightarrow \T_{\Eis}$.
\label{lemma:upper2}
\end{lemma}

\begin{Proof}
Our proof is along the same lines as Lemma~\ref{lemma:upper}.
Suppose that $\rho: \Gal(\Qbar/\Q) \rightarrow \GL_2(\F_p[x]/x^2)$
is upper triangular, and let $\G/\Ot_{K}$ be the associated
finite flat group scheme. Let $\psi$ be the character corresponding
to the upper left hand corner of the representation. Then the
Galois representation $\psi|_{K}$ gives rise to a sub-representation
of the generic fibre of $\G$, and thus to a finite flat subgroup scheme
$\Ht$ of $\G$. The generic fibre of the group scheme $\Ht$ has a filtration by two
copies of the Galois
module $\F_p(\omega^k)$. Since $\F_p(\omega^k)$ extends to
a unique group scheme $\Ht'$ over $\Ot_{K}$
by Corollary~\ref{cor:unique}, the scheme $\Ht$
is an element of $\Ext^1(\Ht',\Ht')$. Moreover, $\Ht$
admits generic fibre descent data to $\Q_p$. It follows from
Corollary~\ref{cor:theoremY}
 that all such extensions split over the maximal
unramified extension  $K^{ur}$. In particular, we see that the
character $\psi: \Gal(\Qbar/\Q) \rightarrow \F_p[x]/x^2$ must be  
the $k$th power of the cyclotomic character.  
Thus $\rho$ must have the shape:
$$\left(\begin{matrix} \chi^k & \st \\ 0 & 
\chi^{1-k} \end{matrix}\right) \in \GL_2(\F_p[x]/x^2).$$
The kernel of $\rho$ therefore defines a $(\Z/p \Z)^2$
extension of $\Q(\zeta_p)$ on which $\Gal(\Q(\zeta_p)/\Q)$
acts via $\chi^{2k-1}$. Yet since the 
 $\chi^{2k-1}$ eigenspace of the class group
of $\Q(\zeta_p)$ is trivial (by assumption),
the maximal $\chi^{2k-1}$ extension of $\Q(\zeta_p)$
has order $p$. Thus
we have a contradiction, and no upper non-trivial
upper triangular
deformation of $\rhobar$ to $\GL_2(\F_p[x]/x^2)$ exists. 
That $R$ is generated by
traces follows as in Lemma~\ref{lemma:upper}, and that
there exists a surjective map $R \rightarrow \T_{\Eis}$ follows
similarly as in Lemma~\ref{lemma:surj}. It suffices to
prove that the representations $\rho_f$ do actually
come from (inverse limits of)  finite flat group schemes over  $K$.
This follows from~\cite{Gross}, Corollary 12.5, since
the associated abelian varieties  acquire good reduction 
over the extension (there denoted by) $M$ of ramification
degree $e | p+1$.
$\qed$
\end{Proof}

\begin{lemma} There is no non-trivial element of
$\D(\Z/p^2 \Z)$ such that $\rho$ is upper triangular.
\label{lemma:lower2}
\end{lemma}

\begin{Proof}
Consider such a representation $\rho$. As in
Lemma~\ref{lemma:upper2}, consider the character $\psi|_K$,
and the (uniquely) associated finite flat group scheme $\Ht/\Ot_K$.
By considering the filtration on the generic fibre of $\Ht$,
we once more infer that $\Ht \in \Ext^1(\Ht',\Ht')$. Since the
generic fibre of $\Ht$ is not killed by $p$, it follows
that $\Ht$ itself not killed by $p$. But by Corollary~\ref{cor:theoremZ},
there are no extensions at all of  $\Ht'$ by $\Ht'$ not killed by $p$!
Thus we are done. $\qed$
\end{Proof}

\medskip

We conclude from Lemma~\ref{lemma:upper2} and
Lemma~\ref{lemma:lower2} that the ideal of reducibility
$I$ is maximal, and thus that $R$ is a discrete valuation
ring, and $R \simeq \T_{\Eis}$.
To complete the
proof of 
Theorem~\ref{theorem:maingross}, we must
prove that $\T_{\Eis} \simeq \Z_p$ when $p \equiv 3 \mod 4$
and $k = (3p-1)/4$. 
Since $R$ is a $\Z_p$-algebra of characteristic zero, it
suffices to prove that $R$ does not admit a surjective
map to $\F_p[x]/x^2$, or equivalently that $\D(\F_p[x]/x^2)$
is empty. We have already seen that $\D(\F_p[x]/x^2)$
does not contain any upper triangular elements. Suppose that
it contains an irreducible representation $\rho$.
Let $\varrho$ be the twist of $\rho$ by $\chi^{-k}$.
Then  $\varrhobar$ is the representation:
$$\left(\begin{matrix} 1 & \st \\ 0 &
\omega^{(p-1)/2} \end{matrix}\right).$$
Thus the kernel of $\varrhobar$ is a degree $p$
extension of $F = \Q(\sqrt{-p})$, the quadratic
subfield of $\Q(\zeta_{p})$.
Let $L = F.\Q(\zeta_p)$ be the kernel of $\varrhobar$,
and $H$ the kernel of $\varrho$.
There is an exact sequence:
$$0 \rightarrow \Gal(H/L) \rightarrow \Gal(H/\Q)
\rightarrow \Gal(L/\Q) \rightarrow 0.$$
We claim this sequence is a semidirect product.
The element of order $2$  in $\Gal(L/\Q)$ lifts
uniquely. Since $\Gal(L/\Q) \subset \GL_2(\F_p[x]/x^2)$,
we see that the $2$-Sylow subgroup acts an involution.
Any lifting of the order $p$ element
of $\Gal(H/\Q)$ (which will have order $p$) therefore
provides a splitting. 
By assumption $\varrho|_L$ is not upper triangular.
It follows that $\Gal(H/L)$ has order $p^3$,
and moreover $H$ must have a subfield
$E$ such that $\Gal(E/F)$ has order $p^2$ and
$\Gal(F/\Q)$ acts on $\Gal(E/F)$ as $-1$.
Yet this is a contradiction, since $F$ admits
at most (in fact exactly) one extension of degree $p$
of this form. Thus $\varrho$ and $\rho$ must be
upper triangular, a contradiction. 

\medskip

In general we have not ruled
out the possibility that $\T_{\Eis}$ is always
$\Z_p$, but one suspects that this is a feature
of the limited range of computation available.

\section{Breuil Modules}
\label{section:breuil}

Throughout this section we shall freely
refer to the results and notation of~\cite{Breuil}.
A reference for Breuil modules killed by $p$
 as $k[u]/u^{ep}$-modules
is~\cite{BCDT}, and we also use some theorems
of Savitt~\cite{Savitt} to determine certain extensions killed
by $p$ of
finite flat group schemes  that admit
generic fibre descent data to $\Q_p$.
Our general approach will be 
to prove some statements
about extensions of finite flat group scheme over an
arbitrary tamely ramified discrete valuation ring. The techniques and
technology are essentially due to Breuil~\cite{Breuil},
following results of Fontaine. 
Let $K/\Q_p$ be a tamely ramified extension (of
arbitrary degree).

\begin{theorem} Let $\G$ be a finite flat group scheme
of order $p^2$ over $\Ot_K$.
Then $\G$ sits inside an exact sequence
$$0 \rightarrow \G_{s,b} \rightarrow \G \rightarrow \G_{r,a}
\rightarrow 0$$
of Oort--Tate group schemes.
If $\G$ is not killed by $p$ then there
exists a non-trivial morphism of group schemes
$\G_{r,a} \rightarrow \G_{s,b}$. Moreover, if
$\G$ is not killed by $p$ then
 $\G[p]$ is finite flat if and only if
$\G$ is \'{e}tale or multiplicative, or equivalently if and only
if
$\G$ isomorphic to $\Z/p^2 \Z$ or $\mu_{p^2}$ over $\Ot_{K^{ur}}$.
For any pairs $(r,a)$ and $(s,b)$ for which there does
exist a non-trivial map $\G_{r,a} \rightarrow \G_{s,b}$ 
that is not an isomorphism, there exists
a corresponding extension $\G$ not killed by $p$ if and only
if $r \ge ps$ or $(e-s) \ge p(e-r)$.
\label{theorem:thZ}
\end{theorem}

This is a combination of Lemma~\ref{lemma:conrad},
Corollary~\ref{cor:main} and Lemma~\ref{lemma:examplegeneral}.
As an application, we have the following:

\begin{cor} \label{cor:theoremZ} Let $K$
be a tamely ramified extension of degree $e = p +1$.
$\G$ be a finite
flat group scheme of order $p^2$ such that the generic
fibre is an extension of $\F_p(\omega^k)$ by $\F_p(\omega^k)$.
Suppose moreover that $k \not\equiv 0,1,(p-1)/2,(p+1)/2 \mod p-1$.
Then $\G$ is killed by $p$. In particular, the generic fibre
is killed by $p$.
\end{cor}

\begin{Proof} From the classification of group schemes of
order $p$ (see for example Example 5.2
of \cite{BCDT}, Theorem~\ref{theorem:rank1} and
Corollary~\ref{cor:unique})
we find that for the $k$ outside the exceptional
listed set, the Galois representation $\F_p(\omega^k)$ arises
from a unique finite flat group scheme of order $p$. It follows
from Theorem~\ref{theorem:thZ}  that either $\G = \G[p]$ or $\G$ is \etale or
multiplicative.  Since $\F_p(\omega^k)$ is not \etale or multiplicative,
the result follows. $\qed$
\end{Proof}

\subsection{Definitions}

We rely extensively on the reference~\cite{Breuil}.
Let $p$ be an odd prime, and let $e$ be an integer coprime to $p$.
Let $\F \subset \Fbar_p$, $W = W(\F)$, $K_0 = W \otimes \Q_p$.
Let $K$ be a totally
tamely ramified extension of $K_0$. 
Let $\pi$ be a uniformizer of $\Ot_K$ with minimal
polynomial
$E(u)=u^{e} + p$. Let $v_n = v(p^n!)$.
Let $S$ be the $p$-adic completion of
$$W[u,X_n], \ \text{where} \ X_n = \frac{u^{e p^n}}{p^{v_n}} \ 
\text{for}  \ n \ge 1.$$
Let $\Fil^1 S$ be the
$W$-submodule of $S$ topologically
generated by $Y_n = E(u)^{p^n}/p^{v_n}$ for all $n$.
There is an isomorphism:
$$S/\Fil^1 S \simeq \Ot_K, \qquad u \mapsto \pi.$$
Let $\phi$ be the unique additive map $S \rightarrow S$,
semilinear with respect to the absolute Frobenius on $W$, 
continuous for the $p$-adic topology, compatible with
the divided powers, and satisfying $\phi(u) = u^p$. Let
$\phi_1 = \frac{\phi}{p} |_{\Fil^1 S}$, and
let $S_n = S/p^n$.

\begin{df} \label{df:one} The category of Breuil modules
$($denoted by
$\ '\ModS)$ consists of triples $(\M,\Fil^1 \M,\phi_1)$
such that \begin{itemize}
\item  $\M$ is an $S$-module
\item  $\Fil^1 \M$ is a
sub $S$-module of $\M$ containing
$\Fil^1 S \cdot \M$,
\item $\phi_1$ is a $\phi$-semilinear map
$\Fil^1 \M \rightarrow
\M$ such that for
all $s \in \Fil^1 S$ and $x \in \M$,
$\phi_1(sx) = \phi_1(s) \varphi(x)$,
where $\varphi(x) = \phi_1(E(u) x)/\phi_1(E(u))$.
\end{itemize}
A map between Breuil
modules is a map $\M \rightarrow \M'$ such
that the induced map on  $\Fil^1 \M$ has image in $\Fil^1 \M'$,
and commutes with $\phi_1$.
The category
$\ModSp$ is the category of Breuil modules with
$\M$ a free $S_1 = S/p$-module of finite rank
such that $\phi_1(\M)$ generates $\M$ as an $S$-module.
\end{df}

Note that in our case, $E(u) = u^e + p$ and
$\phi_1(E(u)) = 1 + u^{ep}/p = 1 + X_1$ (which is 
a unit in $S$). If $\M$ is killed by $p^n$ it
is still important that
one take
$s \in \Fil^1 S$ in this last condition
rather than $s \in \Fil^1 S \cdot S/p^n$.
This is because
$\phi_1(s)\mod p$ does not depend only
on $s\mod p$. For example, ``$\phi_1(u^e)$'' $  = X_1 \ne \phi_1(u^e + p) \mod p$.

\medskip We use the notation $\varphi$ instead of 
$\phi$ used in~\cite{Breuil} to avoid any confusion with
the map $\phi$ defined on $S$.

\begin{theorem} Let $\ModS$ be the sub-category of
$\ '\ModS$ generated from extensions by  $\ModSp$. Then
$\ModS$ is anti-equivalent to the category of finite
flat group schemes over $\Ot_K$.
\end{theorem}

\medskip

Note that $S_1 \simeq \F[u,X_n]/(u^{ep},X^p_n)$.
Moreover, $\Fil^1 S_1 = \Fil^1 S \cdot S_1 = (u^e,X_n) S_1$,
so 
$$S_1/\Fil^1 S_1 \simeq \Ot_K/p \simeq \F[u]/u^{e}.$$

\subsection{Rank One Breuil Modules}

Let $\A(r,a)$ be the Breuil module corresponding
to the following data:
$$\A(r,a) = S_1 \ve,  \quad  \Fil^1 \A(r,a) = (u^r,X_n) \ve, 
\quad
\phi_1(u^r \ve)
= a \ve.$$
\begin{lemma} Any rank one Breuil module killed by
$p$ is isomorphic to $\A(r,a)$ for some $r \le e$
and $a \in \F$. Moreover, there exists a non-trivial map
$\A(s,b) \rightarrow \A(r,a)$ if and only
if $s \equiv r \mod p -1$ and $a/b \in \F^{\times (p-1)}$,
and an isomorphism if and only
if $s = r$ and $a/b \in \F^{\times (p-1)}$.
\end{lemma}

\begin{Proof}
This is essentially~\cite{BCDT}
Example~5.2, and~\cite{Breuil}
Prop.~2.1.2.2,
except for the claim that $\phi(u^r \ve_1)$
can be chosen to equal $a \ve_1$ rather than 
$a x \ve_1$ for some unit $x \in S_1$ satisfying
$x \equiv 1 \mod (X_n)$.
Changing variables by $\ve' = y \ve$, it suffices
to solve the equation
$$\frac{y}{\phi(y)} = x.$$
Since $x \equiv 1 \mod (X_n)$, $\phi^{(n)}(x) = 1$ for
sufficiently large $n$. Thus one can choose
$$y = \prod_{n=0}^{\infty} \phi^{(n)}(x).$$ 
$\qed$
\end{Proof}

\medskip

Group schemes of order $p$ are also classified by
their Oort--Tate parameters \cite{OT}. The following lemma
records that these parameters are (essentially) $(r,a)$.
\begin{lemma} \label{lemma:oort}
The Breuil module $\A(r,a)$
corresponds to the Oort--Tate group scheme
$\G_{r,a}$ over $\Ot_K$
where the affine algebra of $\G_{k,a}$  is equal to
$$\Ot_K[X]/(X + \pi^{e-k} \widetilde{a})$$
where $\widetilde{a}$ is a lift of $a$ to $W(\F)$.
\end{lemma}

\begin{Proof} See for example~\cite{Breuil}
Prop.~3.1.2.
\end{Proof}

\subsection{An Example} 
\label{section:example}

Recall that the category
$\ModFIS$ is the subcategory of
$\ModS$ consisting of Breuil modules such that
$\M \simeq \oplus S_{n_i}$. The category
$\ModFIS$ corresponds to finite flat group schemes
$\G$ such that $\G[p^i]$ is finite flat for all $i$.
We now construct an explicit example of a Breuil
module in $\ModS$ that does not lie in $\ModFIS$.
This example is in~\cite{Breuil}, but we feel that
it serves as a useful example of the modules that
shall be considered in the next section.
Let $e = p-1$.
Recall that $\A(e,1)$, $\A(0,1)$ are the
rank one Breuil modules given by the following data:
$$\begin{array}{lll}
\A(e,1) = S_1 \ve_1, & \Fil^1 \A(e,1) = \Fil^1 S
\cdot \ve_1 = (u^e,X_i) \ve_1, & \phi(u^e \ve_1) = \ve_1 \\
\A(0,1) = S_1 \ove_2,   & \Fil^1 \A(0,1) = S_1 \ove_2,  & \phi_1(\ove_2) = \ove_2
\end{array}$$
In addition, let $\B(0,1)$ be the Breuil module given by
$$\B(0,1) = S_2 \ve_2,  \quad \Fil^1 \B(0,1) = S_1 \ve_2,   \quad \phi_1(\ve_2)
= \ve_2.$$
We use the notation $\ove_2$ for a generator of $\A(0,1)$
to highlight the fact that $\A(0,1)$ is a quotient
of $\B(0,1)$ where the quotient map (multiplication by $p$)
sends $\ve_2$ to $\ove_2$.
The Breuil module $\A(e,1)$ corresponds to the finite
flat constant group scheme $\Z/p\Z$,  whilst the module 
$\A(0,1)$ corresponds to $\mu_p$
(See Lemma~\ref{lemma:oort} above for an identification
of rank one Breuil modules with Oort--Tate group schemes).
Define $\psi: \A(0,1) \rightarrow \A(e,1)$ by
$\psi(\ove_2) = u^p \ve_1 = u^{e+1} \ve_1.$
We check that
$$\phi_1(\psi(\ove_2)) = 
\phi_1(u^p \ve_1) = 
\phi_1(u \cdot u^{e} \ve_1) = \phi(u) \ve_1 = u^p \ve_1
 = \psi(\ove_2) = \psi(\phi_1(\ove_2)).$$
The Breuil module $\B(0,1)$
is an extension of $\A(0,1)$ by  $\A(0,1)$
(corresponding to $\mu_{p^2}$), and there
is a natural map $\iota:\A(0,1) \rightarrow \B(0,1)$ given by
$\ve_2 \rightarrow \ove_2$.  Consider the map
$$\iota +  \psi: \A(0,1) \rightarrow \B(0,1) \oplus \A(e,1),
\qquad \ve_2 \mapsto p \ve_2 -  u^p \ve_1,$$
and let $\Lt$ be the cokernel.
Abstractly it is the quotient
of $S_2 \oplus S_1$ by the element
$(p,-u^p)$.
There is an injective map $\A(e,1) \rightarrow \Lt$ given by
$\ve_1 \mapsto \ve_1$, which extends to a map of
Breuil modules.
The quotient of $\Lt$ by $\A(e,1)$
is the quotient of $S_2 \oplus S_1$ by $\ve_1$
and $p \ve_2 - u^p \ve_1$. Since together these elements
 generate
the module $((p \ve_2,0),(0,\ve_1))$, the quotient is
$\B(0,1)/p \simeq \A(0,1)$.
Thus we have an exact sequence of Breuil modules:
$$0 \rightarrow \A(e,1) \rightarrow \Lt \rightarrow \A(0,1) \rightarrow 0.$$
This extension of Breuil modules corresponds to an
exact sequence of group schemes (note the anti-equivalence):
$$0 \rightarrow  \mu_p \rightarrow \G \rightarrow \Z/p \Z 
\rightarrow 0.$$

\subsection{Extension Classes}

Suppose that $\Ht_1$ and $\Ht_2$  are finite flat group schemes
over $\Ot_K$ of order $p$. 
Let $\G$ be an extension of $\Ht_1$ by $\Ht_2$, and 
let $\M$ be the associated Breuil module. 

\begin{lemma} $\M$ is generated by $($at most$)$ two elements
as an $S$-module. $\M$ is a quotient of $S_2 \oplus S_1$.
\end{lemma}

\begin{Proof} Let $\A_i$ be the Breuil module associated to
$\Ht_i$. There is an exact sequence
$$0 \rightarrow \A_1 \rightarrow \M \rightarrow \A_2 \ra 0.$$
Call the quotient map $\psi$. 
Suppose that $\A_1$  and $\A_2$
are  generated by $\ve_1$ and $\overline{\ve}_2$.
Then
$\M$ is generated by $\ve_1$ and $\ve_2:= \psi^{-1} 
\overline{\ve}_2$.
Moreover, $p \ve_1 = 0$, and $p \ve_2 \in S_1 \ve_1$ so
$p^2 \ve_2 = 0$. $\qed$
\end{Proof}

We conclude that
$$\M = (\ve_1 S_1 \oplus \ve_2 S_2)/I$$
for some $S$-submodule $I$.
\begin{lemma} $I$ is generated by $p \ve_2 - \eta \ve_1$
for some $\eta \in S_1$.
\end{lemma}

\begin{Proof}
The element $p \ve_2$ lies in the kernel
of $\M \rightarrow \A_2$. Thus $p \ve_2$
lies in the image of $\A_1$, and
thus
$p \ve_2 = \eta \ve_1$ for some $\eta \in S_1$.
Suppose that $\beta \ve_2 = \alpha \ve_1$.
If $\beta = p \gamma$, then
$$\gamma (p \ve_2 - \eta \ve_1)
= \beta \ve_2 - \gamma \eta \ve_1 = (\beta \ve_2 - \alpha \ve_1)
+ (\alpha - \gamma \eta) \ve_1$$
and so $(\alpha - \gamma \eta) \ve_1 = 0$ in
$\M$.
Yet $\A_1$ injects into $\M$, so $\gamma \eta = \alpha$,
and 
$$(\beta \ve_2 - \alpha \ve_1) = \gamma (p \ve_2 - \eta \ve_1).$$
Thus we may assume that $p \nmid \beta$.
As an abstract abelian group,
 $S_2 \simeq (\Z/p^2 \Z)^{\infty}$, and thus the image
of $\beta$ is non-zero in $S_1$. The
map from $\M \rightarrow \A_2$ sends $\ve_1$ to $0$.
Thus the image $\overline{\ve}_2$
of $\ve_2$ is  killed by
$\overline{\beta}$,
which is a contradiction, since $\M$ is surjective and
$\A_2$ is free. $\qed$
\end{Proof}

Thus we have an isomorphism of $S$-modules
$$\M =  (S_1 \oplus \ve_2 S_2)/(p \ve_2 - \eta \ve_1).$$

\medskip

Recall that the category
$\ModFIS$ consists of Breuil modules such that
$\M \simeq \oplus S_{n_i}$. The category 
$\ModFIS$ corresponds to finite flat group schemes
$\G$ such that $\G[p^i]$ is finite flat for all $i$.
We conclude that $\G = \G[p]$ if and only if $\eta = 0$,
and that if $\G[p] \ne \G$, then $\G[p]$ 
 is finite flat if and only
if $\eta$ is a unit in $S$.

\medskip

Let us now choose $\A_1 = \A(r,a)$ and $\A_2 = \A(s,b)$.
where $r$ and $s$ are integers $\le e$.
There is an induced exact sequence
$$0 \rightarrow (u^r,X_n) \A(r,a) \rightarrow
\Fil^1 \M \rightarrow (u^s,X_n) \A(s,b) \rightarrow 0.$$
Recall that $\Fil^1 S_2$ is generated by $Y_n = E(u)^{p^n}/p^{v_n}$
for all $n$, 
where $v_n = v(p^n!)$. 
\begin{lemma} Suppose that $p > 2$. Let $X_0 = u^e$,
and let $X_n = u^{ep^n}/p^{v_n}$.
Then
$$Y_n \equiv
 X_n  + p \prod_{i=0}^{n-1} X^{p-1}_{i} \mod p^2 S.$$
In the module $\M$, $X_n \ve_2 \in \Fil^1 \M$ for
all $n \ge 1$.
\end{lemma}

\begin{Proof} The congruence follows by induction
from the identity
$$\frac{(a + bp)^p}{p} = \frac{a^p}{p} + a^{p-1} b p  \mod p^2$$
and the fact that $Y_{n+1} = Y^p_{n}/p$.
Since $e(p-1) \ge e \ge r$, and since $u^r \ve_1 \in \Fil^1 \M$,
it follows that 
$$p X^{p-1}_0 \ve_2  = u^{e(p-1)} \eta \ve_1 \in \M.$$
Thus as $Y_n \ve_2$ is automatically in $\Fil^1 \M$
the inclusion $X_n \ve_2 \in \Fil^1 \M$ follows from the
congruence by induction.
$\qed$ 
\end{Proof}

\medskip

\begin{lemma} If $\M$ is an extension of $\A(s,b)$ by $\A(r,a)$ then
$\M = S_1 \oplus S_2/(p \ve_2 - \eta \ve_1)$. Moreover,
$$\Fil^1 \M = (u^r \ve_1,u^s \ve_2 + x \ve_1,X_n \ve_1,
X_n \ve_2), \ \ i \ge 1,$$
and $\phi_1$ is defined as follows:
$$\phi_1(u^r \ve_1) = a \ve_1, \qquad
\phi_1(u^s \ve_2 + x \ve_1) = b \ve_2.$$
\end{lemma}

\begin{Proof}
We clearly have that
$\Fil^1 \M \cap \A(r,a) = (u^r,X_n) \ve_1$, and
$u^s \ve_2 + x \ve_1 \in \Fil^1 \M$ for some $x \in S$.
Consider a general element $\gamma$ of $\Fil^1 S$.  Since
$X_n \ve_2 \in \Fil^1 \M$, we may assume after subtracting
some element of $(X_n) \ve_2$ that 
$\gamma = \alpha u^s \ve_2 + \beta \ve_1$. Then
$$\alpha u^s \ve_2 + \beta \ve_1 - \alpha (u^s \ve_2 + x \ve_1)
= (\beta - \alpha x) \ve_1 \in \Fil^1 \M.$$
Thus $(\beta - \alpha x) \ve_1 = \Fil^1 \A(r,a)$, and $\gamma$
is in the span of $u^s \ve_2 + x \ve_1$, $(X_n) \ve_2$, and
$\Fil^1 \A(r,a)$. Hence
$$\Fil^1 \M = (u^r \ve_1,u^s \ve_2 + x \ve_1,X_n \ve_1, 
X_n \ve_2), \ \ i \ge 1.$$
Moreover,
$$\phi_1(u^r \ve_1) = a \ve_1, \phi_1(u^s \ve_2 + x \ve_1) = b \ve_2 + z 
\ve_1.$$
Replacing $\ve_2$ by $\ve_2 + y \ve_1$ for suitable
$y$ we may assume
(having
changed $x$ and $\eta$ appropriately) that $z = 0$.
This proves the lemma. $\qed$
\end{Proof}

\medskip

We also note that 
that $\phi_1(\Fil^1 \M)$ as defined above generates $\M$ 
as an $S$-module. However, we have not shown
that $\phi_1$ is well defined, and so we have
not yet constructed a Breuil module. In fact, the obstructions
to defining $\phi_1$ will severely limit the possible
extension classes.

\medskip

Fix once and for all an element $y = \eta - u^{e-s} x$.

\medskip
Let us compute $\varphi(\ve_2)$ and $\varphi(\ve_1)$. 
Recalling that $y = \eta - u^{e-s} x$ we find that
$$(u^e + p) \ve_2 = u^e \ve_2 + \eta \ve_1 = u^{e-s}(u^s \ve_2
+ x \ve_1) - u^{e-s} x \ve_1 + \eta \ve_1
= u^{e-s}(u^s \ve_2
+ x \ve_1) + y \ve_1.$$
It follows that $y \ve_1 \in \Fil^1 \A(r,a)$.
We find that
$$\phi_1(E(u) \ve_2) = b u^{p(e-s)}  \ve_2 + \phi_1(y \ve_1).$$
and so
$$\varphi(\ve_2) =
\frac{\phi_1(E(u) \ve_2)}{\phi_1(E(u))}
=  \frac{b u^{p(e-s)} \ve_2  + \phi_1(y \ve_1)}{1 + X_1}.$$
A similar computation shows that
$$\varphi(\ve_1) = \frac{\phi_1(E(u) \ve_1)}{\phi_1(E(u))}
= \frac{\phi_1(u^{e-r} u^r \ve_1)}{1 + u^{ep}/p}
= \frac{a \cdot  u^{p(e-r)} \ve_1}{1 + X_1}.$$

When $\M = \Lt$, we see that
$\Fil^1 \Lt$ is generated by $\ve_2$ and $u^r  \ve_1$.
Moreover, $\phi_1(\ve_2) = \ve_2$ and $\phi_1(u^r \ve_1) = \ve_1$.
Thus $y = \eta = u^p$, and $x = 0$. Moreover,
$s = 0$ and $b = 1$
and so
$$\varphi(\ve_2) = \frac{u^{pe} \ve_2 + u^p \ve_1}{1 + X_1}
= \frac{p X_1 \ve_2 + u^p \ve_1}{1 + X_1}
= \frac{X_1 \eta \ve_1 + u^p \ve_1}{1 + X_1} = 
\frac{(1 + X_1) u^p \ve_1}{1 + X_1} = u^p \ve_1.$$

\medskip 

Returning now to the general situation, we shall derive
some relations between $\eta$, $x$ and $y$ by computing
$\phi_1(u^s Y_n \ve_2 + x Y_n  \ve_1)$ in two different ways. 
On the one hand we have that
$$\phi_1(u^s Y_n \ve_2) = \phi(u^s) \phi_1(Y_n) \varphi(\ve_2)
= u^{ps} Y_{n+1} \cdot \frac{b u^{p(e-s)} \ve_2  + 
\phi_1(y \ve_1)}{1 + X_1}$$ 
and 
$$\phi_1(x Y_n \ve_1) = \phi(x) \phi_1(Y_n) \varphi(\ve_1)
= \phi(x) Y_{n+1} \cdot  \frac{a u^{p(e-r)} \ve_1}{1 + X_1}.$$
On the other hand, $u^s Y_n \ve_2 + x Y_n \ve_1
= Y_n(u^s \ve_2 + x \ve_1)$, and so
$$\phi_1(Y_n(u^s \ve_2 + x \ve_1)) = \phi(Y_n) b \ve_2.$$
Now $\phi(Y_n) = Y^{p}_n = p Y_{n+1}$, and $p \ve_2 =
\eta \ve_1$. Thus we find that
$$Y_{n+1} b \eta \ve_1 = Y_{n+1} \left(\frac{b u^{pe} \ve_2 + u^{ps} \phi_1(y \ve_1)
+ a \phi(x) u^{p(e-r)} \ve_1}{1 + X_1} \right).$$
We make two simplifications. First, $u^{ep} \ve_2 = X_1 p \ve_2 = X_1 \eta \ve_1$.
Thus both sides are multiples of $\ve_1$, and we may replace
$Y_{n+1}$ by $X_{n+1}$.
Second, the annihilator of $X_{n}$ in $S_1$ is $X^{p-1}_n$. Thus the
annihilator of $(X_n)_{n=k}^{m}$ is $\prod_{n=k}^{m} X^{p-1}_n$
and the annihilator of $(X_n)_{n=k}^{\infty}$ is trivial. 
Thus if $X_{n+1} \alpha = X_{n+1} \beta$ for all sufficiently
large $n$, $\alpha = \beta$. Applying this to our formula,
and 
multiplying through by $(1 + X_1)$ we find that
$$ b \eta (1 + X_1) \ve_1 =  X_1 b \eta \ve_1  + u^{ps} \phi_1(y \ve_1)
+ a \phi(x) u^{p(e-r)} \ve_1$$
or
$$b \eta \ve_1 = u^{ps} \phi_1(y \ve_1)
+ a \phi(x) u^{p(e-r)} \ve_1.$$
Thus we have proven:
\begin{theorem} \label{theorem:rest} If $\M$ is well defined as a Breuil
module then
$$b \eta \ve_1 = u^{ps} \phi_1(y \ve_1)
+ a \phi(x) u^{p(e-r)} \ve_1.$$
Moreover $\eta$ is divisible by $(u^p)^{\min\{s,e-r\}}$.
\end{theorem}

Suppose that $\M = \Lt$. Then
$x = 0$, $e = p-1$, and $y = \eta = u^p$. Thus
the theorem is consistent with the identity
$u^p = \eta = \phi_1(\eta \ve_1)$.

\medskip

\begin{lemma} \label{lemma:nose}  
There is an inclusion
$\eta \in \F[u]/u^{ep}.$
\end{lemma}

\begin{Proof} We divide our proof into two
cases. First we consider the case where $(r,s) \ne (e,0)$.
Since $\phi(X_n) = X^{p}_n = 0$, it
is clear that $\phi(x) \in  \F[u]/u^{ep}$. 
Moreover, $\phi_1(X_i X_j \ve_1) = \phi(X_i) \phi_1(X_j \ve_1)
= X^p_n \phi(X_j \ve_1) = 0$ for the same reason. Thus the only terms that contribute
to coefficients of $\eta$ that do not lie in
$\F[u]/u^{ep}$ are of the form $u^m X_n$.
Let $m$ be the infimum (minimum) over all $n \ge 1$ such
that the coefficient of $u^m X_n$ in $y$ is non-zero.
The corresponding coefficient of $\eta$ is
$u^{p(e - r + s + m)} X_{n+1}$, and thus if $m \ge r-s$ we
are done. Suppose otherwise. 
Then the minimum $m$ over all $n \ge 1$ such that the
coefficient of $u^m  X_{n}$ in $\eta$ is non-zero
is $p(e + s - r + m)$. The minimum $m$
over all $n \ge 1$ such that the
coefficient of $u^m X_n$ in $u^{e-s} x$ is non-zero is
trivially at least 
at least $e-s$. Since 
$y = \eta - u^{e-s} x$ we conclude that
$$m \ge \min\{p(e + s - r + m),e-s\}.$$
Since $m < r-s \le e-s$, it must be the first inequality
that is satisfied. Equivalently,
$$(1-p) m \ge p(e + s - r).$$
Since $e \ge r$, the RHS is non-negative unless $r = e$
and $s = 0$. On the other hand, the LHS is negative 
unless  $m \ge 0$. Thus either we are done or
$(r,s) = (e,0)$ and $m=0$.  Let us now assume we are
in that case.
There is an identity
$y = \eta - u^e x$. Since $y \ve_1 \in \Fil^1 \A(e,a)$ and
$u^e \in \Fil^1 S_1$, it follows that $\eta \in \Fil^1 S_1$.
Now $\phi_1(y \ve_1) = \phi_1(\eta \ve_1) - a \phi(x)$  and so
$$ b \eta \ve_1 = \phi_1(y \ve_1) + a \phi(x)
= \phi_1(\eta \ve_1).$$ 
As above, since $\phi_1(X_i X_j \ve_1) = \phi(X_i) \phi(X_j \ve_1)
= 0$, the only terms contributing to $\eta$ that
do not lie in $\F[u]/u^{ep}$ are coefficients
of $y$ of the form
$u^m X_n$ with $n \ge 1$. Let $n$ be the smallest integer such
that $u^m X_n$ is a non-zero coefficient of $\eta$. Then since
$$\phi_1(u^i X_{j}) = u^{ip} \phi_1(Y_{j+1}) \varphi(\ve_1)
= \frac{X_{j+1} u^{ip}}{1 + X_1}$$
we see that (since $j +1 >  m$) that
$\phi_1(\eta)$ does not have any coefficients of
the form $u^n X_m$, a contradiction. Thus $\eta \in \F[u]/u^{ep}$.
$\qed$
\end{Proof}

\medskip

Write 
$$x = \sum_{k=0}^{ep-1} \alpha_k u^k \mod (X_n),
\qquad
y = \sum_{k=0}^{ep-1} \beta_k u^k \mod (X_n).$$
Then (noting  by Lemma~\ref{lemma:nose} that $\eta \in \F[u]/u^{ep}$) the equality
$\eta = y + u^{e-s} x$ becomes
$$\sum_{k=0}^{ep-1} \gamma_k u^k = \eta = \sum_{k=0}^{ep-1} u^k(\beta_k + 
\gamma_{k+s-e}).$$
Since $y \ve_1 \in \Fil^1 \A(r,a)$, we must have
$\beta_k = 0$ for $k < r$. 
Applying the equality
of Theorem~\ref{theorem:rest} we find that
$$b  \eta = \sum_{k=0}^{ep-1} a u^{pk}(\phi(\beta_{k+r-s}) + \phi(\alpha_{k+r-e}))
= a \sum_{k=0}^{ep-1} u^{pk}\phi(\beta_{k+r-s} + \alpha_{k+r-e}),$$
where $\phi$ is Frobenius on $\F$.
The two expressions for $\eta$ lead to the following relations:
\begin{itemize}
\item $\gamma_k = \beta_k + \alpha_{k+s-e}$
\item $\gamma_k = 0$ if $p \nmid k$.
\item $b \gamma_{pk} = a \phi(\beta_{k+r-s} + \alpha_{k+r-e})$.
\end{itemize}
In particular we see that
$b \gamma_{pk} = a \phi(\gamma_{k + r-s})$.

\begin{lemma} if $\M$ is a Breuil module, then
$\eta = 0$ unless $(r-s) = k(p-1)$ for some
$k \ge 0$
and $a/b \in \F^{\times (p-1)}$.
If $r = s$, then $\eta = 0$
unless $s = 0$ or $r = e$. If $\eta$ 
is non-zero, then up to
an element of $\F^{\times}$,
$\eta = u^{kp}$. Finally, there is
an inequality $k \ge \min\{s,e-r\}.$
\label{lemma:vital}
\end{lemma}

\begin{Proof}
First note that $\gamma_{pk}$
is non-zero if and only if $\gamma_{k + (r-s)}$
is non-zero, since Frobenius is injective in $\F$.
Let $k$ be the smallest integer
such that $\gamma_{pk} \ne 0$. Then $k \ge 0$
and $k$ is the smallest integer such
that 
$\gamma_{k + (r-s)} \ne 0$. It follows 
that $k + (r-s) = pk$, and so $(r-s) = k(p-1)$.
The equality $b \gamma_{pk} = a \phi(\gamma_{pk}) \ne 0$
implies that $a/b = \phi(c)/c \in \F^{\times (p-1)}$.
One finds (by considering the \emph{second} smallest
$k$ such that $\gamma_{pk} \ne 0$) 
that no other coefficients of $\eta$ are non-zero,
and thus
$\eta$ is a multiple of $u^{pk}$.
 If $r=s$ then Theorem~\ref{theorem:rest}
implies that $\eta$ is divisible by $u$ (unless
$s = 0$ or $r = e$),
but the $k$ satisfying $(r-s) = k(p-1)$ is $k = 0$,
thus $\eta = 0$. The final inequality follows
from Theorem~\ref{theorem:rest}. $\qed$
\end{Proof}

\begin{cor} \label{cor:main}  Suppose that $\G$ is an extension:
$$0 \rightarrow \G_{s,b} \rightarrow \G \rightarrow \G_{r,a} 
\rightarrow 0.$$
 of Oort--Tate group schemes that is \emph{not} killed  by $p$.
Then there is a non-trivial morphism $\G_{r,a} \rightarrow \G_{s,b}$.
\end{cor}

\begin{Proof}
Suppose there existed such an exact sequence. Then there must
exist an exact sequence of Breuil modules:
$$0 \rightarrow \A(r,a) \rightarrow \M \rightarrow \A(s,b)
\rightarrow 0.$$
If $\G$ is not killed by $p$, then $\M$ is not killed by $p$,
and thus $\eta \ne 0$. By Lemma~\ref{lemma:vital} then the
restrictions on $(r,s)$ and $(a,b)$ are exactly the requirements
that there exist a map $\A(s,b) \rightarrow \A(r,a)$
(see~\cite{BCDT}, Lemma 5.2.1). Moreover, this map
cannot be an isomorphism unless $r=s$, which 
(from Lemma~\ref{lemma:vital}) implies that $r=s=0$ or
$r=s=e$. In this case $\G_{r,a}$ and $\G_{s,b}$ are
either both multiplicative or both \etalestop, which
implies that $\G$ is either multiplicative or \etalestop.
$\qed$ \end{Proof}

\medskip

if $r$ and $s$ are both $0$ or both $e$ there exist extensions
not killed by $p$: one can take
$\mu_{p^2}$ and $\Z/p^2 \Z$ respectively.

\begin{lemma} Suppose that $r$ and $s$ are
integers $\le e$ such that
$(r-s) = k(p-1)$ and $k \ge \min\{s,e-r\}$. 
Let $a$ and $b$ be elements of $\F$ such that $(b/a) = (c)^{p-1}$.
Then there exists a non-trivial extension of the form:
$$0 \rightarrow \G_{s,b} \rightarrow \G \rightarrow \G_{r,a}
\rightarrow 0.$$
\label{lemma:examplegeneral}
\end{lemma}
\begin{Proof} Suppose that $k \ge s$. Then one may
explicitly let
$\M = S_1 \oplus S_2/(p \ve_2 - c  u^{pk} \ve_1)$. Since
 $kp = r + k-s \ge r$, it follows that
$u^{kp} \ve_1$ is a multiple of $u^{r} \ve_1$. Let
$$\Fil^1 \M = (u^r \ve_1, u^s \ve_2, X_n \ve_1, X_n \ve_2),$$
and finally define  $\phi_1$ as follows:
$$\phi_1(u^r \ve_1) = a \ve_1, \qquad
\phi_1(u^s \ve_2) = b \ve_2.$$
If $k \ge e-r$ then we still define 
$\M = S_1 \oplus S_2/(p \ve_2 - c  u^{pk} \ve_1)$. Now that
$k \ge e-r$ we define $\Fil^1 \M$ as follows
$$\Fil^1 \M = (u^r \ve_1, u^s \ve_2 + c  u^{k-e+r}
 \ve_1, X_n \ve_1, X_n \ve_2)$$
where
$$\phi_1(u^r \ve_1) = a \ve_1, \qquad
\phi_1(u^s \ve_2 + c  u^{k-e+r} \ve_1) = b \ve_2.$$
Note that $u^{e-s}(u^s \ve_2 + c  u^{k-e+r} \ve_1)
= u^e \ve_2 + c  u^{pk} \ve_1 = (u^e + p) \ve_2$.
One verifies in both cases that $\M$ defines a Breuil
module.
$\qed$
\end{Proof}

\subsection{Finite Flat Group Schemes Killed by $p$}

If one restricts to finite flat group schemes killed
by $p$, the theory of Breuil modules can be
significantly simplified. In particular, instead
of working with $S_1$-modules, it suffices
to work with $S_1/(X_n) = k[u]/u^{ep}$ modules,
and replace $\M$ by $\M \otimes_{S_1} k[u]/u^{ep}$.

\medskip

What Breuil modules $\A(r,a)$ admit generic descent
data to $\Q_p$? The answer for a general
tamely ramified extension is provided in~\cite{Savitt}.
We restrict the statement to the case of interest,
namely when
 $K/\Q_p$ is a tamely ramified
Galois extension of $\Q_p$ with $e =p+1$.
We have the following:

\begin{theorem} \label{theorem:rank1} The Breuil module $\A(r,a)$ admits
generic fibre descent data to $\Q_p$ if and only if
$2$ divides $r$ and $a \in \F^{\times}_p$.
Let $\xi_a$ be the unramified character
of $\Gal(\Qbar_p/\Q_p)$ given by
$\xi(\sigma) = \sigma a^{1/(p-1)}/a^{1/(p-1)}$.
Then
the associated Galois representation on the
(descended) generic
fibre is given by $\xi_a \F_p(\omega^k)$, where $r
\equiv 2 - 2k \mod (p-1)$.
\end{theorem}

\begin{Proof} This follows from the
calculations in
in~\cite{Savitt}, in particular
Definition 5.1, Proposition 5.3, and Theorem 6.3. 
Note in our
setting and Savitt's notation we have
$U = -1$, $V = -1$, and $x' = x/(e,p-1) = 2$. $\qed$
\end{Proof}

\begin{cor} \label{cor:unique}
If $k \ne 0,1,(p-1)/2,(p+1)/2$, then
there is a unique finite flat group scheme of
order $p$ with generic fibre given by
$\F_p(\omega^k)$.
\end{cor}

\begin{Proof}  Any such finite flat group scheme
obviously has generic fibre descent data. 
We see that such a representation
forces $a$ to be $1$, and $r \not\equiv 0,2 \mod (p-1)$.
Yet since $r$ is even and less than $p+1$,
this determines $r$ exactly. $\qed$
\end{Proof}

\begin{lemma} \label{lemma:theoremX}
 Let $\rhobar$ be as in
section~\ref{section:gross}. Then
$\rhobar$ uniquely determines a finite
flat group scheme $\G/\Ot_K$. 
\end{lemma}

\begin{Proof}
Since $\F_p(\omega^k)$ and $\F_p(\omega^{1-k})$ correspond to
unique finite flat group schemes, the lemma follows
from a standard application of the $5$-Lemma
(see~\cite{BCDT}, in particular the proof of
Lemma 4.1.2). $\qed$
\end{Proof}

Now we turn to extensions of Breuil modules killed
by $p$.  An easy computation (\cite{BCDT},
\cite{Savitt}) shows that
the extensions of $\A(r,a)$ by itself are classified
by an element 
$h \in u^{\max(0,2r-e)} k[u]/u^{r+1}$. In general,
these will correspond to finite flat group schemes
whose generic fibre does not descend to $\Q_p$.
However, we have the following:

\begin{theorem} \label{theorem:savitt} The
space of extensions $\Ext^1(\A(r,a),\A(r,a))$
killed by $p$ has dimension $1$ over $\F_p$.
\end{theorem}

\begin{Proof} It follows from Theorem~7.5 of~\cite{Savitt}
that $h$ can be taken to have degree less than $r$.
Moreover, (since $k_1 = k_2$ in the notation of~\cite{Savitt})
the only possible non-zero term is $u^r$ (which is zero)
or the constant term which lies in $\F_p$. $\qed$
\end{Proof}

\begin{cor} \label{cor:theoremY}  Suppose that $K$ is
tamely ramified of degree $p+1$, that $k \ne 0,1,(p-1)/2,(p+1)/2$.
Let $\Ht/\Ot_K$ be a finite flat group scheme over $\Ot_K$ with
generic fibre $\F_p(\omega^k)$. Then the only extensions of $\Ht$ by
itself which admit generic fibre descent data to $\Q_p$ become
unramified over some finite unramified extension of $K$.
\end{cor}

\begin{Proof}
The space of extensions is one dimensional by theorem~\ref{theorem:savitt}.
Thus it suffices to observe that the Galois module $\F_p(\omega^k)$
(over $\Q_p$) 
admits an extension by itself that splits over the degree
$p$ unramified extension of $\Q_p$. Thus if $F$ denotes this
unramified extension, then by faithfully flat descent from
$\Ot_{F.K}$ to $\Ot_K$ we obtain a non-trivial
(in fact, $p-1$ non-trivial)  extensions of
$\Ht$ by $\Ht$ in the category of finite flat group schemes
with generic fibre descent data, which splits over some
unramified extension of $K$, and thus we are done.
$\qed$
\end{Proof}


\begin{thebibliography}{99}

\bibitem{Chenevier}
J.~Bella\"{\i}che, G.~Chenevier.
\emph{Lisset\'{e} de la courbe de Hecke
de $\GL_2$ aux points Eisenstein critiques},
Preprint.

\bibitem{Breuil}
C.~Breuil.
\emph{Groups $p$-divisibles, groupes finis
et modules filtr\'{e}s},
Ann. of Math. (2) \bf152\rm (2000), no. 2, 489--549. 

\bibitem{BCDT}
C.~Breuil, B.~Conrad, F.~Diamond, R~.Taylor.
\emph{On the modularity of elliptic curves over $\Q$: wild $3$-adic exercises},
J. Amer. Math. Soc. \bf14\rm (2001), no. 4, 843--939.

\bibitem{CE}
F.~Calegari, M.~Emerton.
\emph{Ramification of Hecke Algebras at 
Eisenstein Primes},
to appear in Invent. Math.

\bibitem{calstein}
F.~Calegari, W.~Stein.
\emph{Conjectures about discriminants of Hecke Algebras},
to appear in the proceedings of ANTS VI.

\bibitem{Conrad}
B.~Conrad.
\emph{Ramified Deformation Problems},
Duke Math. J. \bf97\rm (1999), no. 3, 439--513.



\bibitem{CDT}
B.~Conrad, F.~Diamond, R.~Taylor.
\emph{Modularity of certain potentially Barsotti-Tate Galois representations},
J. Amer. Math, Soc. \bf12\rm (1999), no. 2, 521--567.

\bibitem{Flach}
M.~Flach,
\emph{A finiteness theorem for the symmetric square of an elliptic curve},
Invent. Math. \bf109\rm (1992), no. 2, 307--327.


\bibitem{FL}
J.~Fontaine, G.~Laffaille.
\emph{Construction de repr\'{e}sentations $p$-adiques},
Ann. Sci. \'{E}cole Norm. Sup. (4) \bf15\rm (1982), no. 4, 547--608 (1983)

\bibitem{Greenberg}
R.~Greenberg.
\emph{Iwasawa Theory --- Past and Present},
Class field theory---its centenary and prospect (Tokyo, 1998), 335--385.

\bibitem{Gross}
D.~Gross, J.~Lubin.
\emph{The Eisenstein Descent on $J_0(N)$},
Invent. Math. \bf83\rm (1986), no. 2, 303--319.


\bibitem{Lang}
S.~Lang.
\emph{Cyclotomic Fields},
GTM \bf59\rm, Springer Velag, 1978, 253 pp.

\bibitem{eisenstein}
B.~Mazur.
\emph{Modular curves and the Eisenstein ideal},
Publ.~Math.~IHES \bf47\rm (1977), 33--186

\bibitem{OT}
F.~Oort, J.~Tate.
\emph{Group schemes of prime order},
Ann. Sci. \'{E}cole Norm. Sup. (4) \bf3\rm (1970) 1--21.

\bibitem{Raynaud}
M.~Raynaud.
\emph{Sch\'{e}mas en groupes de type $(p,\dots, p)$},
Bull. Soc. Math. France \bf102\rm (1974), 241--280.

\bibitem{Ribet}
K.~Ribet.
\emph{A modular construction of 
unramified $p$-extensions of $\Q(\zeta_{p})$},
Invent. Math. \bf 34\rm, (1976), no. 3, 151--162.

\bibitem{Papier}
K.~Ribet, E.~Papier.
\emph{Eisenstein ideals and $\lambda$-adic representations},
J. Fac. Sci. Univ. Tokyo Sect. IA Math. \bf28\rm 
(1981), no. 3, 651--665 (1982).



\bibitem{Savitt}
D.~Savitt.
\emph{Modularily of some potentially Barsotti--Tate
representations},
Compos. Math. \bf140\rm (2004), no. 1, 31--63.

\bibitem{SW}
C.~Skinner, A.~Wiles.
\emph{Ordinary representations and modular forms},
Proc. Nat. Acad. Sci. U.S.A. \bf94\rm (1997), no. 20, 10520--10527.

\end{thebibliography}
\end{document}